\documentstyle{amsppt}

\magnification=\magstep1
\baselineskip=22pt
\parindent=8mm
\TagsOnRight
\font\fontA=cmcsc10 
 
%
%=========================================================================
%
\pageno=1
\document
\vskip 2cm

\pageno=1

\centerline{\bf $C^*$-ALGEBRAS WHOSE EVERY $C^*$-SUBALGEBRA IS AF}

\medskip
\medskip
\medskip 
\medskip 

\rm
%\centerline{\fontA MASAHARU KUSUDA}
\centerline{MASAHARU KUSUDA}

\medskip 
\medskip 

\centerline  {\fontA Department of Mathematics,} \par
\centerline  {\fontA Faculty of Engineering Science,} \par
\centerline  {\fontA Kansai University,} \par
\centerline  {\fontA Yamate-cho 3-3-35, Suita, Osaka 564-8680, Japan.}\par
\centerline{  e-mail address: kusuda\@ipcku.kansai-u.ac.jp}

\vskip1cm

\topmatter
\abstract 
Let $A$ be a $C^*$-algebra. 
It is shown that the following conditions are equinvalent:

\noindent (1) $A$ is scattered, 

\noindent (2) every $C^*$-subalgebra of $A$ is AF,  

\noindent (3) every $C^*$-subalgebra of $A$ has real rank zero. 

\endabstract
\endtopmatter
\medskip

\medskip
\medskip
\medskip
\medskip

\comment
\topmatter
\endtopmatter
\endcomment

\medskip 
\medskip 
\medskip 

\medskip 

\noindent {\it Key words and key phrases}: scattered $C^*$-algebra, 
AF $C^*$-algebra, $C^*$-algebra of real rank zero
\medskip

\newpage 

\head{1. Introduction}\endhead

An AF $C^*$-algebra was introduced by Bratteli, which is defined 
as the inductive limit of a sequence of finite-dimensional 
$C^*$-algebras. 
As is well-known, every hereditary $C^*$-subalgebra of 
an AF $C^*$-algebra is also AF. But in general, 
an AF $C^*$-algebra contains $C^*$-subalgebras which are not AF. 
For example, the $C^*$-algebra $C(2^\omega)$ of all 
continuous functions on the Cantor set $2^\omega$ is AF because 
the Cantor set $2^\omega$ is a totally disconnected topological space. 
But $C(2^\omega)$ contains a $C^*$-subalgebra which 
is isomorphic to the $C^*$-algebra $C(I)$  of all continuous functions 
on the closed interval $ I = [0, 1]$, and  $C(I)$ 
is not AF.  In \S 2, 
by an AF $C^*$-algebra we will mean 
a $C^*$-algebra $A$ in which 
any finite set of elements can be 
approximated arbitrarily closely in norm by elements of a finite-dimensional 
  $C^*$-subalgebra of $A$. This definition  
has a slightly wider sense 
than the usual one introduced by Bratteli, 
but an AF $C^*$-algebra in our sense is recently often  adopted 
as the definition of an AF $C^*$-algebra. 
Then 
any scattered $C^*$-algebra is AF in our sense above (\cite{10, Lemma 5.1}). 
However, every AF $C^*$-algebra is not scattered. 
See \cite{9, Theorem 2.2} for conditions for  AF $C^*$-algebras to become 
scattered. 
The purpose of this paper is to characterize 
$C^*$-algebras whose every $C^*$-subalgebra  is AF. 
In fact, in \S2, we will show that a $C^*$-algebra $A$ is scattered 
if and only if every $C^*$-subalgebra of $A$ is AF. 

Here recall that 
a topological space is called {\it scattered} (or dispersed) 
if every non-empty subset 
necessarily contains an isolated point.
For a compact Hausdorff space $\Omega$, it is shown in \cite{12} that 
$\Omega$ is scattered if and only 
if every Radon measure on $\Omega$ is atomic. 
As a non-commutative generalization of a scattered 
compact Hausdorff space,  the notion of a scattered $C^*$-algebra  
was introduced independently by Jensen \cite{7} and Rothwell \cite{13}. 
We say that a $C^*$-algebra $A$ is {\it scattered}  if every positive 
linear functional on $A$ is the countable sum of 
pure positive linear functionals 
on $A$, or equivalently, $A$ is of type I and the spectrum 
$\widehat A$ of $A$ is a scattered topological space 
equipped with the Jacobson topology (\cite{8, Corollary 3}, 
\cite{13, Theorem 3.8}). 
The reader is referred to  \cite {3}, \cite{6}, \cite{7}, \cite{8}, 
\cite{9}, \cite{13}  
for other equivalent conditions on scattered $C^*$-algebras.  

By the way, AF $C^*$-algebras of Bratteli's sense have real rank zero 
which  was defined  by Brown-Pedersen \cite{5}, and 
there are $C^*$-algebras with real rank zero which are not AF. 
In \S2, 
we will show also that every $C^*$-subalgebra of $A$ is AF in our sense 
if and only if every $C^*$-subalgebra of $A$ has real rank zero. 
Thus we will see that scattered $C^*$-algebras can be completely 
characterized by the condition that every $C^*$-subalgebra should have 
the same property to be  AF or to have real rank zero. 
This leads us to the question of what kind of $C^*$-algebra can be 
characterized by some  property the $C^*$-subalgebras possess. 
In \S3, we will show that the class of 
type I $C^*$-algebras is  one of answers for 
such a question. In fact, we show that 
a $C^*$-algebra $A$ is of type I if and only if 
every $C^*$-subalgebra of $A$ is strongly amenable, 
if and only if every $C^*$-subalgebra of $A$ is nuclear.

\medskip

\head{2. Approximately finite dimensionality and real rank zero}\endhead

In this section, the term AF $C^*$-algebra has a slightly wider sense 
than the usual one which 
is assumed to be separable and to have identity. 
We will ease those  restrictions. In fact, 
we say that a  $C^*$-algebra $A$ is  
an {\it approximately finite-dimensional} $C^*$-algebra 
(or simply {\it AF}) if 
for any elements $x_1, x_2, \cdots, x_n \in A$ and any $\varepsilon > 0$, 
there exists a finite-dimensional $C^*$-subalgebra $B$ such that 
$\Vert x_k - y_k\Vert  < \varepsilon $ with some $y_k \in B$ 
for all $k = 1, 2, \cdots, n$. 
In this definition, any scattered $C^*$-algebra is AF 
(\cite{10, Lemma 5.1}).

For a compact Hausdorff space $\Omega$, 
we denote by 
$C(\Omega)$ the $C^*$-algebra of all continuous 
functions on $\Omega$. 
Throughout this section, we denote by $2^\omega$ the Cantor set. 
Note that the Cantor set $2^\omega$ is totally disconnected and 
is not scattered. 
Hence $C(2^\omega)$ is AF, but it is not a scattered $C^*$-algebra. 

Here we recall that a topological space $X$ is called 0-{\it dimensional} 
if each point of $X$ has a neighborhood base consisting 
of open-closed sets. Equivalently, $X$ is 0-dimensional 
if and only if for each point $x \in X$ and closed  set $F$ not 
containing $x$, there is an open and closed set containing 
$x$ and not meeting $F$. 
It is well-known that 
a locally compact Hausdorff space is 0-dimensional if and only if 
it is totally disconnected.

For a $C^*$-algebra $A$, we always denote  by  
$A_1$ the $C^*$-algebra obtained from $A$ 
by adjunction of an identity. 
We say that a $C^*$-algebra $A$ has {\it real rank zero} if and only if 
every self-adjoint element in $A_1$ can be approximated by 
an invertible self-adjoint element in $A_1$, or equivalently, 
$A$ has the (FS) property , that is, 
the set of all self-adjoint elements with finite spectrum 
of $A$ is dense in the set of all self-adjoint elements of $A$ 
(\cite{5, Theorem 2.6}).

\medskip
\medskip

\proclaim{Lemma  2.1} Let $\Omega$ be a  compact Hausdorff space 
and let $C(\Omega)$ be the $C^*$-algebra of all continuous functions 
on $\Omega$. 
If every $C^*$-subalgebra of $C(\Omega)$ has real rank zero, 
then $C(\Omega)$ is scattered. 
\endproclaim

\medskip
\noindent {\it Proof}.  
Assume that $C(\Omega)$ is not a scattered $C^*$-algebra. 
Then $\Omega$ is not scattered as a topological space 
(see \cite{8, Corollary 3}). 
On the other hand,  
since $C(\Omega)$ has real rank zero by assumption, 
it follows from \cite{5, 1.1} that 
$\Omega$ is 0-dimensional. 
Hence it follows from \cite{12, 2. Main Theorem} that 
there exsists a continuous map $\pi$ from $\Omega$ onto 
the Cantor set $2^\omega$. 
Hence 
we obtain that
$C(2^\omega)$ can be  embedded into $C(\Omega)$  as a $C^*$-algebra 
via the induced map $\widetilde \pi$ defined by 
$\widetilde \pi(f) = f \circ \pi$ for $f \in C(2^\omega)$. 
Since there is a continuous map from 
$2^\omega$ onto the closed interval $I = [0, 1]$, 
$C(I)$ can be embedded into  $C(2^\omega)$ as a $C^*$-algebra. 
Hence $C(I)$ can be embedded into  $C(\Omega)$ as a $C^*$-algebra. 
Since $I$ is not totally disconnected, it is not 0-dimensional. 
Hence $C(I)$ does not have real rank zero. 
This is a contadiction because every $C^*$-subalgebra of $C(\Omega)$ 
has real rank zero. Thus we conclude that 
$C(\Omega)$ is a scattered $C^*$-algebra. \qed

\medskip
\medskip
 Let $\Omega$ be a compact Hausdorff space. We remark that 
if $C(\Omega)$ is scattered, then 
it has real rank zero. But the converse is false in general.  
A counterexample is $C(2^\omega)$, which is seen from the proof of 
Lemma 2.1. 
On the other hand, from Lemma 2.1 above and this remark , we see 
that every $C^*$-subalgebra of $C(\Omega)$ has real rank zero if 
and only if  every $C^*$-subalgebra of $C(\Omega)$ is scattered.

\medskip
\medskip
For a $C^*$-algebra $A$, 
we denote by $\widehat A$  
the spectrum of $A$ which is the set of 
all  equivalence classes of nonzero irreducible 
representations equipped with the Jacobson topology.  
As is well known, 
$\widehat A$ is locally compact, but it is not necessarily 
a Hausdorff space. If $A$ is unital, 
$\widehat A$ is compact. 

Note that every $C^*$-subalgebra of 
a scatrered $C^*$-algebra $A$ is scattered 
(\cite{6, Theorem} or \cite{3, Theorem 1 and Theorem 3}).

\proclaim{Lemma  2.2} A  $C^*$-algebra $A$ is scattered 
if and only if 
every separable abelian $C^*$-subalgebra of $A$ is scattered. 
\endproclaim

\medskip
\noindent {\it Proof}.  Suppose that every separable 
abelian $C^*$-subalgebra of $A$ is scattered. 
If $A$ is not unital, we consider the $C^*$-algebra $A_1$ obtained by 
adding an identity 1 to $A$. If $A_1$  is scattered, so is $A$. 
Without loss of generality, thus we may assume that $A$ is unital. 
A local characterization of a scattered $C^*$-algebra 
 is that each self-adjoint element $h$ of the $C^*$-algebra has a 
countable spectrum $\text{Sp}_A(h)$ of $h$ in the $C^*$-algebra 
(\cite{6, Theorem}). Now we use this characterization. 

Let $h$ be any self-adjoint element in $A$. 
Since the $C^*$-subalgebra $C^*(h)$ generated by $h$ and $1$ 
is separable and abelian, by assumption $C^*(h)$ is scattered. 
Since $C^*(h)$ is separable, 
with \cite{7, Theorem 3.1} and \cite{8, Theorem 2} combined, 
we see that the spectrum $\widehat {C^*(h)}$ of $C^*(h)$ is countable. 
Since $\widehat {C^*(h)}$ is homeomorphic to 
$\text{Sp}_A(h)$, $\text{Sp}_A(h)$ is countable. 
Hence it follow from \cite{6, Theorem} that 
$A$ is scattered. 
\qed

\medskip
\medskip

Now  we are in a position to give the main theorem in this section. 
\medskip
\medskip

\proclaim{Theorem 2.3} Let $A$ be a  $C^*$-algebra. 
Then the following conditions are eqivalent.

$(1)$ $A$ is scattered. 

$(2)$ Every $C^*$-subalgebra of $A$ is AF.  

$(3)$  Every separable abelian $C^*$-subalgebra of $A$ is AF.

$(4)$ Every $C^*$-subalgebra of $A$ has real rank zero.

$(5)$  Every separable abelian $C^*$-subalgebra of $A$ 
has real rank zero.

\endproclaim

\medskip
\noindent {\it Proof}. $(1) \Longrightarrow (2)$. 
If $A$ is scattered, then 
every $C^*$-subalgebra  of $A$ is scattered. 
Since a scattered $C^*$-algebra is AF (\cite{10, Lemma 5.1}),  
every $C^*$-subalgebra of $A$ is AF. 

$(2) \Longrightarrow (4)$. It is well known that 
every AF $C^*$-algebra of the sense of Bratelli has real rank zero 
(\cite{5, 3.1}). 
Every AF $C^*$-algebra in the our sense also 
has real rank zero. 
In fact, let $B$ be an  AF $C^*$-algebra in the our sense. 
For any self-adjoint element $h \in B$ and any 
$\varepsilon > 0$, 
there exists an element $a$ in some finite-dimensional 
$C^*$-subalgebra $B_0$ of $B$ such that $\Vert h - a \Vert < \varepsilon$. 
Then we have 
$\Vert h - \frac{a + a^*}{2} \Vert < \varepsilon$. 
Hence $\Vert h - k \Vert < \varepsilon$ for 
some self-adjoint element $k \in B_0$.
Since $B_0$ is a finite-dimensional $C^*$-algebra, 
$k$ has finite spectrum. 
Thus the set of all self-adjoint elements with finite spectrum 
of $B$ is dense in the set of all self-adjoint elements of $B$. 
This shows that $B$ has real rank zero (\cite{5, Theorem 2.6}), 
from which condition (4) follows.

$(4) \Longrightarrow (5)$. This is obvious.

$(5) \Longrightarrow (1)$. 
We assume that every separable abelian $C^*$-subalgebra of $A$ 
has real rank zero. 
Take any separable abelian $C^*$-subalgebra $B$ of $A$. 
Then by assumption, $B$ has real rank zero. 
First we assume that $B$ has  identity. 
Since every $C^*$-subalgebra $C$ of $B$ is 
also a separable abelian $C^*$-subalgebra of $A$, 
$C$ has real rank zero. 
It hence follows from Lemma 2.1 that $B$ is scattered. 
Then $A$ is scattered by Lemma 2.2. 

Next we assume that $B$ does not have identity. 
 Let $B_1$ be the $C^*$-algebra obtained by adjoining an identity 1 to $B$. 
Since $B$ has real rank zero, 
$B_1$ also does so by the definition of real rank. 
From the discussion of the unital $C^*$-algebra case above, 
we see that $B_1$ is scattered. Hence $B$ is scattered. 
 Then $A$ is scattered by Lemma 2.2.

$(2) \Longrightarrow (3)$. This is obvious. 

$(3) \Longrightarrow (5)$. 
The proof is similar to that of $(2) \Longrightarrow (4)$. 
Or this implication easily follows also from 
\cite{5, 3.1}. 
\qed

\medskip

\head{3. Amenability}\endhead

\medskip

In this section, 
we consider the question of what kind of $C^*$-algebra 
can be characterized by some property all the $C^*$-subalgebras  possess, 
besides scattered $C^*$-algebras we have considered in the previous section.  

First we briefly review  the definition of amenability of 
a $C^*$-algebra. Let $A$ be a $C^*$-algebra. 
We say that a Banach space $X$ is called a Banach $A$-module if 
it is a two-sided $A$-module and there exists  a constant $K > 0$ 
such that for any $a \in A$ and $x \in X$ we have 
$$ \Vert ax \Vert \leqq K\Vert a \Vert\Vert x\Vert \quad
\text{and} \quad \Vert xa \Vert \leqq K\Vert x \Vert\Vert a\Vert. $$

For a  Banach $A$-module $X$, the dual space $X^*$ of $X$ becomes a 
Banach $A$-module by 
$$ (af)(x) = f(xa) \quad
\text{and} \quad (fa)(x) = f(ax) $$
for all $a \in A$ and $f \in X^*$. 
A bounded linear map $D$ from $A$ into $X^*$ 
is called a derivation if it satisfies 
$D(ab) = aD(b) + D(a)b$ for $a, b \in A$. 
If a derivation $D$ from $A$ into $X^*$ 
is given by  $D(a) = af - fa$  for some $f \in X^*$, 
then $D$ is called inner.   
A $C^*$-algebra $A$ is called {\it amenable} if every bounded derivation 
from $A$ into $X^*$ is inner for all Banach $A$-modules $X$. 
When $A$ does not have identity, we  denote again by $A_1$ 
the $C^*$-algebra obtained from $A$ 
by adjunction of an identity $1$. 
Then we can make a Banach $A$-module $X$ into an $A_1$-module 
by $x1 = 1x = x$ for $x \in X$, and $D$ can be extended to 
a derivation form $A_1$ into $X^*$ by defining $D(1) = 0$.

A $C^*$-algebra $A$ is called {\it strongly amenable} 
if for every  Banach $A$-module $X$, every bounded derivation $D$ 
from $A$ into $X^*$ is given by 
$D(a) = af - fa$ with some $f$ 
in the weakly${}^*$ closed convex hull of 
$\{ -D(u)u^* | u \in U(A_1) \}$, where $U(A_1)$ 
denotes the unitary group of $A_1$. 
Obviously a strongly amenable $C^*$-algebra is amenable. 

Recall that a $C^*$-algebra $A$ is said to be {\it nuclear} if 
$A \otimes_{\max} B = A \otimes_{\min} B$ 
for any $C^*$-algebra $B$, 
where $\otimes_{\max}$ and $\otimes_{\min}$ 
denote the maximal $C^*$-tensor product and 
the minimal $C^*$-tensor product, respectively. 
Equivalently, 
a $C^*$-algebra $A$ is  nuclear if 
the identity map on $A$ is nuclear, that is, 
for any elements $x_1, x_2, \cdots, x_k \in A$ and any $\varepsilon > 0$, 
there exist a natural number  $n$ and 
completely  positive contractions $\Phi : A \to M_n(\Bbb C)$ and 
$\Psi : M_n(\Bbb C) \to A$ 
such that $\Vert \Psi \circ \Phi (x_i) - x_i \Vert < \varepsilon$ 
for all $i = 1, 2, \cdots, k$, where $M_n(\Bbb C)$ 
is the $n \times n$-matrix algebra (e.g., \cite{1, IV.3.1.5}). 
It is well-known that 
a $C^*$-algebra is amenable if and only if it is nuclear 
(\cite{1, IV.3.3.15}).

In \cite{2, Theorem 1}, Blackadar showed that 
a type I $C^*$-algebra contains a non-nuclear $C^*$-subalgebra. 
This result yields the implication $(4) \Longrightarrow (1)$ in the Theorem 3.1 below. 
To show \cite{2, Theorem 1}, 
he constructed an example of a non-nuclear $C^*$-subalgebra 
in the type I $C^*$-algebra. The construction of the example is not easy. 
In the below, we present an alternative proof of 
the implication $(4) \Longrightarrow (1)$ in the Theorem 3.1 
based on a well-known theorem of Glimm \cite{11, 6.7.4}. 
\medskip

\proclaim{Theorem 3.1}  Let $A$ be a $C^*$-algebra. 
Then the following conditions $(1) - (5)$ are equivalent.

%\noindent 
$(1)$ $A$ is of type I. 

%\noindent 
$(2)$ Every $C^*$-subalgebra of $A$ is strongly amenable. 

%\noindent 
$(3)$ Every separable $C^*$-subalgebra of $A$ is strongly amenable.

%\noindent 
$(4)$ Every $C^*$-subalgebra of $A$ is nuclear.

%\noindent 
$(5)$ Every separable $C^*$-subalgebra of $A$ is nuclear. 
\endproclaim

\medskip
\noindent {\it Proof}. 
$(1) \Longrightarrow (2).$ 
This follows from the fact that every type I $C^*$-algebra is 
strongly amenable (\cite{4, Theorem 7.9}). 

$(2) \Longrightarrow (4).$ 
 This follows from the fact that every strongly amenable $C^*$-algebra is 
amenable, hence nuclear. 

$(4) \Longrightarrow (1).$ 
Assume that $A$ is not of type I. 
Then there is a $C^*$-subalgebra $B$ of $A$ such that 
there exists a surjective homomorphism $\pi$ from $B$ 
onto the CAR algebra $\Cal F$ (see \cite{11, 6.7.4}). 
Let $C^*_r(F_2)$ be the reduced group $C^*$-algebra of the free group $F_2$ 
on two generators. Since 
$C^*_r(F_2)$ is a separable exact $C^*$-algebra (cf. \cite{14, 2.5.3}), 
it follows from \cite{1, IV.3.4.18 (iv)} that 
there is a $C^*$-subalgebra $C$ of $\Cal F$ such that 
there exists a surjective homomorphism $\rho$ from $C$ 
onto $C^*_r(F_2)$.  Since a homomorphic image 
of a nuclear $C^*$-algebra is nuclear and since 
$C^*_r(F_2)$ is not nuclear (cf. \cite{1, II.9.4.6} or \cite{14, 2.4}), 
$C$ is not a nuclear $C^*$-subalgebra of $\Cal F$. 
Hence $\pi^{-1}(C)$ is not a nuclear $C^*$-subalgebra of $B$. 
Thus we see that $A$ contains the $C^*$-subalgebra $\pi^{-1}(C)$ 
which is not nuclear. This is a contradiction. 

$(2) \Longrightarrow (3) \Longrightarrow (5).$ These are trivial.

$(5) \Longrightarrow (4).$  Let $B$ be a $C^*$-subalgebra of $A$. 
Take any elements $x_1, x_2, \cdots, x_k$ in $B$ 
and any $\varepsilon > 0$. 
Let $B_0$ be the $C^*$-subalgebra of $B$ 
generated by $x_1, x_2, \cdots, x_k$. 
Since $B_0$ is separable, it follows from Condition (5) that 
$B_0$ is nuclear.  Hence 
there exist a natural number $n$ and 
completely positive contractions $\Phi : B_0 \to M_n(\Bbb C)$ and 
$\Psi : M_n(\Bbb C) \to B_0 \subset B$ 
such that $\Vert \Psi \circ \Phi (x_i) - x_i \Vert < \varepsilon$ 
for all $i = 1, 2, \cdots, k$. 
By Arveson's Extension Theorem \cite{1, II.6.9.12}, 
there exists a complete positive contraction $\widehat\Phi : B \to M_n(\Bbb C)$ 
which extends $\Phi$. 
Thus we obtain that 
$\Vert \widehat\Phi \circ \Psi (x_i) - x_i \Vert < \varepsilon$ 
for all $i = 1, 2, \cdots, k$, which shows that
$B$ is nuclear. 
\qed

\enddocument